\magnification= \magstep1
\baselineskip=14pt
\pageno=0
\newcount\knum
\knum=1
\newcount\uknum
\uknum=1
\newcount\znum
\znum=1
\newcount\notenumber

\def\Za{\the\uknum.\the\znum \global\advance\znum by 1}
\def\ukneu{\znum=1\global\advance\uknum by 1}
\font\tbfontt=cmbx10 scaled \magstep0
\font\gross=cmbx10 scaled \magstep2
\font\mittel=cmbx10 scaled \magstep1
\font\ittel=cmr10 scaled \magstep1

\font\boldit=cmbxti10 scaled \magstep0
\font\TT=cmcsc10 scaled \magstep0

\font\sq=cmss10 scaled \magstep0
\font\nirm=cmr9
\font\eightrm=cmr8

\font\eightit=cmti8 scaled \magstep0
\def\sqr#1#2{{\vcenter{\vbox{\hrule height.#2pt\hbox{\vrule width.#2pt
height#1pt \kern#1pt \vrule width.#2pt}\hrule height.#2pt}}}}
\def\square{\mathchoice\sqr34\sqr34\sqr{2.1}9\sqr{1.5}9}
\def\QED{\smallskip\rightline{$\square$ \quad\qquad\null}}

\def\Theorem{\vskip 0.3 true cm\sl {\TT Theorem \Za.}{ }}

\def\Corollary{\vskip 0.3 true cm\sl {\TT Corollary \Za.}{ }}
\def\Lemma{\vskip 0.3 true cm\sl {\TT Lemma \Za.}{ }}
\def\Proof{\TT Proof:{ }\rm}

\def\cz{{\rm C\hskip-4.8pt\vrule height5.8pt\hskip5.3pt}}
\def\kz{{\rm C\hskip-4.0pt\vrule height4.3pt\hskip5.3pt}}  
\def\rz{{\rm I\hskip -2pt R}}
\def\Az{{\rm A\hskip -5.15pt\raise 1pt\hbox{$\scriptscriptstyle
\backslash$}} }
\def\nAz{{\nabla\hskip -6.6pt\raise 2.0pt
\hbox{$\scriptstyle\nabla$}} }
\def\az{{ \hbox{\nirm A}\hskip -4.9pt\raise 0.7pt\hbox{${\scriptscriptstyle
\backslash}$} }}

\def\Ee{{\rm I\hskip -2pt E}}
\def\zz{{\rm Z\hskip -4pt Z}}
\def\Fz{{\rm I\hskip -2pt F}}  
\def\eins{{\rm 1\hskip -2pt I}}

\def\mapright#1{\smash{\mathop{\longrightarrow}\limits^{#1}}}

\def\hookr#1{\smash{\mathop{\hookrightarrow}\limits^{#1}}}

\def\mapdown#1{\Big\downarrow\rlap{$\vcenter{\hbox{$\scriptstyle#1$}}$}}

\def\lim{\mathop{\rm lim}}

\def\End{{\hbox{\rm End}}}

\parindent=0pt
\font\eightrm                = cmr8
\font\eightsl                = cmsl8
\font\eightsy                = cmsy8
\font\eightit                = cmti8

\font\eighti                 = cmmi8
\font\eightbf                = cmbx8
\def\petit{\def\rm{\fam0\eightrm}
\textfont0=\eightrm 
 \textfont1=\eighti 
 \textfont2=\eightsy 
 \def\it{\fam\itfam\eightit}
 \textfont\itfam=\eightit
 \def\sl{\fam\slfam\eightsl}
 \textfont\slfam=\eightsl
 \def\bf{\fam\bffam\eightbf}
 \textfont\bffam=\eightbf 
 \normalbaselineskip=9pt
 \setbox\strutbox=\hbox{\vrule height7pt depth2pt width0pt}
 \normalbaselines\rm}
\newdimen\refindent
\def\begref{\vskip1cm\bgroup\petit
\setbox0=\hbox{[Bi,Sc,So]o }\refindent=\wd0
\let\sl=\rm\let\INS=N}

\def\ref#1{\filbreak\if N\INS\let\INS=Y\vbox{\noindent\tbfontt
References\vskip1cm}\fi\hangindent\refindent
\hangafter=1\noindent\hbox to\refindent{#1\hfil}\ignorespaces}

\long\def\fussnote#1#2{{\baselineskip=9pt
\setbox\strutbox=\hbox{\vrule height 7pt depth 2pt width 0pt}%
\petit\noindent\footnote{\noindent #1}{#2}}}
{\nopagenumbers
\rightline{dg-ga/9601004}
\rightline{January 1996}
\vskip 4cm
\phantom{prelim.version}
\centerline{\gross Supersymmetry and}
\medskip
\centerline{\gross the generalized Lichnerowicz formula}
\vskip 1cm

\centerline{\ittel Thomas Ackermann\footnote{$^*$}{\eightrm
e-mail: ackerm@euler.math.uni-mannheim.de} }
\vskip 0.5cm
\centerline{Wasserwerkstr. 37, D-68309 Mannheim, F.R.G.}
\vskip 1.5cm
\centerline{\vbox{\hsize=5.0 true in \petit\noindent
\bf Abstract.\rm\ A classical result in
differential geometry due to Lichnerowicz [8] is concerned with the
decomposition of the square of Dirac operators
defined by
Clifford connections on a Clifford module
${\cal E}$\ over a Riemannian manifold $M$. Recently, this
formula has been generalized
to arbitrary Dirac operators [2].
In this paper we prove a supersymmetric
version of the generalized Lichnerowicz formula, motivated by
the fact that
there is a one-to-one correspondence between Clifford superconnections
and Dirac operators. We extend this result to obtain a simple
formula for the supercurvature of a generalized Bismut superconnection.
This might be seen as a first step to prove the local index
theorem also for families of arbitrary Dirac operators.}}
\vskip 3.0cm\parindent=1.8cm
\item{Keywords:} \it supersymmetry, generalized
Lichnerowicz formula, generalized Bismut superconnections, local
family index theorem
\par\parindent=0pt
1991 MSC: \it 47F05, 53A50, 58G10\hfill \vfill
\break}\rm
\advance\hsize by -0.5 true in
\advance\vsize by -0.9true in
\advance\hoffset by 0.45 true in
\advance\voffset by 0.7 true in
{\mittel 1. Introduction}
\vskip 0.7cm
A classical result in differential geometry due to Lichnerowicz [8]
is concerned with the decomposition of the square of Dirac operators
defined by Clifford connections on a Clifford module ${\cal E}$\ over a
Riemannian manifold $M$. More precisely, if $M$\ is even dimensional,
$\nabla^{\cal E}\colon \Gamma({\cal E})\rightarrow \Gamma(T^*M\otimes
{\cal E})$\ is such a connection and $D_{\nabla^{\cal E}}:=
c\circ \nabla^{\cal E}$\ denotes the corresponding Dirac operator,
Lichnerowicz's formula states 
$${D_{\nabla^{\cal E}}^2=\triangle^{\nabla^{\cal E}}+{r_M\over 4} +
{\bf c}(R^{{\cal E}/S}_{\nabla^{\cal E}}).}\eqno(1.1)$$
Here $\triangle^{\nabla^{\cal E}}$\
is the connection laplacian associated
to $\nabla^{\cal E}$\ and the endomorphism part is given by
the scalar curvature  
$r_M$\ of $M$\ and 
the image ${{\bf c}(R^{{\cal E}/S}_{
\nabla^{\cal E}})}$\ of the twisting curvature
$R^{{\cal E}/S}_{\nabla^{\cal E}}
\in \Omega^2(M,
\End_{C(M)}({\cal E}))$\ of the Clifford connection $\nabla^{\cal E}$\
with respect to
the quantisation map
${{\bf c}\colon \Lambda^*T^*M\rightarrow C(M)}$, cf. [5].
This formula provides a powerful tool for studying the Dirac operator
$D_{\nabla^{\cal E}}$\ and its relation with the geometry of
the underlying manifold. For example, if $M$\ is a compact spin manifold
with positive scalar curvature, Lichnerowiz used (1.1) to show the vanishing
of the space of harmonic spinors. Moreover, formula (1.1) assumes a
significant r$\hat {\rm o}$le in the proof of the
local Atiyah-Singer
index theorem for such Dirac operators, cf. [6].
\medskip
Recently, Lichnerowicz's formula (1.1) has been extended
to arbitrary Dirac operators: In [2] it was shown that
the square $\widetilde{D}^2$\
of a Dirac operator $\widetilde{D}:=c\circ \widetilde{\nabla}^{\cal E}$\
defined by an arbitrary
connection $\widetilde{\nabla}^{\cal E}\colon \Gamma({\cal E})
\rightarrow \Gamma(T^*M\otimes {\cal E})$\ on the Clifford module
${\cal E}$\ decomposes as
$${\widetilde{D}^2=\triangle^{{\hat\nabla}^{\cal E}}+
{\bf c}(R^{\widetilde{\nabla}^{\cal
E}})+
ev_g\widetilde{\nabla}^{T^*M\otimes {\petit\rm End}{\cal E}}
\varpi_{\widetilde{\nabla}^{\cal E}} +
ev_g(\varpi_{\widetilde{\nabla}^{\cal E}}
\cdot \varpi_{\widetilde{\nabla}^{\cal E}}),}\eqno(1.2)$$
with $\varpi_{\widetilde{\nabla}^{\cal E}}:=-{1\over 2} g_{\nu\kappa}
dx^\nu\otimes
c(dx^\mu)\bigl([\widetilde{\nabla}^{\cal E}_\mu, c(dx^\kappa)] +
c(dx^\sigma)\Gamma^\kappa_{\sigma\mu}\bigr)\in \Omega^1(M,\End {\cal E})$\
and the connection ${\hat\nabla}^{\cal E}:=\widetilde{\nabla}^{\cal E}+
\varpi_{\widetilde{\nabla}^{\cal E}}$.
Here the last two
terms obviously indicate the deviation of the
connection $\widetilde{\nabla}^{\cal
E}$\ being a Clifford connection. Only the second term
in (1.2) is endowed with geometric significance. Of course, if
$\widetilde{\nabla}^{\cal E}$\
is a Clifford connection,
obviously $\varpi_{\widetilde{\nabla}^{\cal E}}=0$\ and therefore
(1.2) reduces to (1.1).
\medskip
Using Quillen's theory of superconnections
on $\zz_2$-graded vectorbundles nowadays it is well-known that 
that any Clifford superconnection
$\Az$\ on $\cal E={\cal E}^+\oplus {\cal E}^-$\ uniquely determines a
Dirac operator ${D_{\Az}}$\
due to the following construction
{$${D_{\az}\colon \Gamma({\cal E})\;\mapright{\az}\;
\Omega^*(M,{\cal E})\;
\mapright{{\petit\bf c}\otimes \eins_{\cal E}}\;\Gamma(C(M)\otimes
{\cal E})\;\mapright{c}\; \Gamma({\cal E}), }\eqno(1.3)$$
i.e. there is a one-to-one correspondence between Clifford
superconnections and}
Dirac operators, see [5]. In contrast to [2] here we emphasize
this approach to Dirac operators. Thus, the first major
purpose of this paper
consists in proving the supersymmetric version
$${D_\az^2=\triangle^{{\hat\nabla}^{\cal E}}+{r_M\over 4}+
{\bf c}\bigl(\Fz(\Az)^{{\cal E}/S}\bigr) +
{\bf c}\bigl({\bar\Az}\bigr)^2 - {\bf c}\bigl(
{\bar\Az}^2\bigr) + ev_g\bigl(\beta({\bar\Az})\cdot
\beta({\bar\Az})\bigr)}\eqno(1.3)$$
of the intrinsic decomposition formula (1.2) in Theorem 4.2.
Here the connection 
${\hat\nabla}^{\cal E}:=\Az_{[1]}+\beta({\bar\Az})$\ is  
determined by the connection part $\Az_{[1]}$\
and the connection-free part ${\bar\Az}\in \bigoplus_{i\ne 1}
\Omega^i(M,\End\;{\cal E})$\ of the Clifford superconnection
$\Az$, $\Fz(\Az)^{{\cal E}/S}$\ denotes the
twisting supercurvature of $\Az$\ and
$\beta({\bar\Az})\in \Omega^1(M,\End\;
{\cal E})$\ is defined by $\beta({\bar\Az}):= dx^k\otimes
{\bf c}\bigl(i(\partial_k){\bar\Az}\bigr)$\ with respect to a
local coordinate frame. Note that Getzler [7] has stated the
generalization of Lichnerowicz formula as   
${D_\az^2\!=\!\triangle^{{\hat\nabla}^{\cal E}}\!+\!{r_M\over 4}\!+\!
{\bf c}\bigl(\Fz(\Az)^{{\cal E}/S}\bigr)\! +\! P(\Az)}$\ without
specifying the endomorphism $P(\Az)\in \Gamma(\End\;{\cal E})$\ in
general. In view of explicit
computations, (1.2) resp. (1.3)
are more convenient, e.g. in applications to physics (cf. [2], [3]).
\medskip
Secondly,
we extend the above formula (1.3) to
a family of Dirac operators $\hbox{\sq D}_\az:=\{\; D_\az\;\vert
\;z\in B\;\}$\ parametrized by a - not necesarily finite
dimensional - manifold $B$. In order to do so we associate to each
family $\hbox{\sq D}_\az$\ a superconnection $\nAz^\az$\ following
Bismut's construcion [4]. In the case of $\hbox{\sq D}_\Az$\
being a family of Dirac operators defined by Clifford connections,
$\nAz^\az$\ reduces to the Bismut superconnection. So we call it
`generalized Bismut superconnection'.
\medskip
The paper is organized as follows: In the next section we establish
our conventions and briefly recall Quillen's superconnection
formalism. In section 3 we establish the canonical projection
${\bf g}\colon \Omega^*(M,{\cal E})\rightarrow
\Omega^1(M,\End\;{\cal E})$\ and show some identities which
are crucial in proving our Theorem 4.2 . Although this
is basic for the theory of Clifford modules we recognize
that there is no presentation available in the
literature. In section 4 we prove our main result, the
supersymmetric version of the generalized Lichnerowicz formula (1.2)
in Theorem 4.2 as we have already mentioned. Finally, in the last
section we extend this formula in Theorem 5.2
to the context of families of
Dirac operators by generalizing Bismut's construction.
\medskip
In another paper [1] we will show how to use Theorem 5.2 to compute
explicitly the Chern character of a generalized Bismut superconnection
which implies the local Atiyah-Singer index theorem
for families of arbitrary Dirac operators.

\vfill\break
{\mittel 2. The superconnection formalism}
\ukneu
\vskip 0.7cm
In this section we briefly introduce some basic results of Quillen
[9] concerning supersymmetry and superconnections. Let $E=E^+\oplus
E^-$\ be a finite dimensional $\zz_2$-graded vector bundle over
a manifold $M$. Then the endomorphism bundle $\End E$\ is canonicaly
graded and - by their total $\zz_2$-gradings - so are the
spaces 
$\Omega^*(M,E)$, $\Omega^*(M,\End E)$\ 
of $E$-vallued and $\End E$-vallued differential forms, respectively.
\smallskip
A superconnection on $E$\ is defined to be an odd parity operator
$\Az\colon \Gamma(E)\!\rightarrow\! \Omega^*(M,E)$\ which satisfies
Leibniz's rule $\Az(fs)=df\otimes s
+ f\Az s$\ for all
$f\in C^\infty(M)$\ and $s\in \Gamma(E)$. Obviously,
$\Az$\ can be extended to an odd operator of $\Omega^*(M,E)$. Furthermore
we can expand $\Az$\ into a sum $\Az=\sum_{i\ge 0}\;\Az_{[i]}$\
of operators $\Az_{[i]}\colon \Omega^*(M,E)\rightarrow \Omega^{*+i}(M,E)$\
such that $\Az_{[1]}$\ is a connection on $E$\ which respects the
grading and $\Az_{[i]}\in \Omega^i(M,\End E)^-$\ for $i\ne 1$. Any
superconnection
$\Az$\ acts on the space $\Omega^*(M,\End E)$\ by $\Az\alpha:=[\Az,
\alpha]$\
for $\alpha\in \Omega^*(M,\End E)$. Here $[\;\cdot\; ,\;\cdot\;]$\
denotes the supercommutator
$${[\alpha,\alpha^\prime]:=\alpha \alpha^\prime - (-1)^{
\vert\alpha\vert\;\vert\alpha^\prime\vert}\alpha^\prime
\alpha}\eqno(2.1)$$
where  $\alpha,\alpha^\prime\in\Omega^*(M,\End E)$. By definition, the
curvature of a superconnection $\Az$\ is $\Fz(\Az):=\Az^2$\
which is an even $C^\infty(M)$-linear endomorphism of $\Omega^*(M,E)$,
i.e. $\Fz(\Az)\in \Omega^*(M,\End E)^+$. Also the supercurvature
$\Fz(\Az)$\ splits into a sum $\Fz(\Az)=\sum_{i\ge 0}\;\Fz(\Az)_{[i]}$\
with $\Fz(\Az)_{[i]}\in \Omega^i(M,\End E)^+$\ given by
$$\eqalign{\Fz(\Az)_{[0]}\; &=\;\Az^2_{[0]}\cr
\noalign{\vskip 0.1cm}   
           \Fz(\Az)_{[1]}\; &=\;[\Az_{[0]},\Az_{[1]}]\cr
           \noalign{\vskip 0.1cm}
           \Fz(\Az)_{[2]}\; &=\;[\Az_{[0]},\Az_{[2]}]+ \Az_{[1]}^2\cr
\noalign{\vskip 0.1cm}
                          &\ \ \vdots \cr\noalign{\vskip 0.1cm}
           \Fz(\Az)_{[i]}\; &=\; [\Az_{[0]},\Az_{[i]}]+[\Az_{[1]},\Az_{[i-1]}]+\;
                         \dots \cr\noalign{\vskip 0.1cm}
                          &\ \ \vdots\cr}\eqno(2.2)$$  
Now let $\nabla^E:=\Az_{[1]}$\ be the connection- and
${{\bar\Az}\in \bigoplus_{i\ne 1}\;\Omega^i(M,\End E)}$\ be the
`connec\-tion-free' part of the superconnection $A$\ such that
the decomposition $\Az=\nabla^E + {\bar\Az}$\ holds.
Because we have $[\nabla^E, \alpha]=d^{\nabla^{End\;E}}\alpha$\ for  
all $\alpha\in \Omega^*(M,\End\;E)$\
where $d^{\nabla^{End\;E}}$\ denotes the induced exterior
covariant derivative
on $\Omega^*(M,\End\;E)$, the above decomposition (2.2) implies 
the following
\Lemma Let
$R^{\nabla^E}\in \Omega^2(M,\End E)$\
denote the curvature of the connection part $\nabla^E:=\Az_{[1]}$\
of a superconnection $\Az\in \Omega^*(M,\End E)$. Then the supercurvature
$\Fz(\Az)\in \Omega^*(M,\End E)^+$\ splits into
$\Fz(\Az)=R^{\nabla^E}+d^{\nabla^{End\; E}} {\bar\Az}+ {\bar\Az}^2$\ with
${\bar\Az}\in \bigoplus_{i\ne 1}\;\Omega^i(M,\End E)$\ the connection-free
part of $\Az$.\rm
\smallskip
The curvature of a superconnection $\Az$\ is thus the sum of the curvature
of the connection $\nabla^E:=\Az_{[1]}$, the exterior covariant
derivative of its connection-free part $\bar \Az$\ with respect
to $\nabla^E$, and ${\bar \Az}^2$.
\smallskip 
Now we turn out attention to a specific class of superconnections on
a Clifford module, namely the Clifford superconnections. Recall, that
a Clifford module 
is a $\zz_2$-graded complex vector bundle
$\cal{E}=\cal{E}^+\oplus \cal{E}^-$\ over
a Riemanniann manifold $M$\ together with a
$\zz_2$-graded left action $c\colon C(M)\times \cal{E} \rightarrow
\cal{E}$\ of the Clifford bundle, i.e. a graded representation
of $C(M)$\fussnote{${^1)}$}{For convenience of
the reader and to fix our conventions we remark that
$C(M)$\ is the bundle of Clifford algebras over $M$\
generated by $T^*M_\kz:=T^*M\otimes_\rz \cz$\ 
with respect to the relations
${v\star w +w\star v=-2 g(v,w)}$\ for
sections $v, w\in \Gamma(T^*M_\kz)$.}. Generalizing the notion
of a Clifford connection (cf. [5]), a Clifford superconnection on
${\cal E}$\ is defined to be
a superconnection $\Az\colon \Gamma({\cal E})^\pm
\rightarrow \Omega^*(M,{\cal E})^\mp$\ which is compatible with
the Clifford action $c$, i.e.
${[\Az, c(a)]=c(\nabla a)}$\
for all $a\in \Gamma(C(M))$.
In this formula, $\nabla$\ denotes the Levi-Civita connection
extended to the Clifford bundle $C(M)$. 
Obviously any Clifford connection defines a Clifford superconnection
with trivial connection-free part ${\bar\Az}=0$. Furthermore
the connection part $\Az_{[1]}$\ of a Clifford
superconnection $\Az$\ determines a Cifford
connection $\nabla^{\cal E}$.
\smallskip
For later use we compare ${\bf c}({\bar \Az}^2)$\ and
${\bf c}({\bar \Az})^2$\ for
the connection-free part $\bar \Az$\ of any
Clifford superconnection $\Az$.
Here ${\bf c}\colon \Omega^*(M, \End\; {\cal E})\rightarrow
\Omega^0(M,\End\;{\cal E})$\ denotes the obvious extension of the
quantisation map ${\bf c}\colon \Lambda^*T^*M\rightarrow C(M)$. Note
that this map ${\bf c}$\ is not a homomorphism of algebras but yields the
identity when
restricted to $\Omega^0(M,\End \;{\cal E})$. Thus, by a simple calculation
we get the following
\Lemma Let $\Az$\ be a Clifford superconnection on a Clifford module
${\cal E}$. Then   
${{\bf c}({\bar\Az})^2
-{\bf c}({\bar\Az}^2)=\sum_{i,j\ge 2}\;\bigl({\bf c}(\Az_{[i]})
{\bf c}(\Az_{[j]}) - {\bf c}(\Az_{[i]}\Az_{[j]})\bigr)}$\ 
where ${\bar \Az}:=\sum_{i\ne 1}\;\Az_{[i]}$\ denotes
the `connection-free' part of $\Az$.
\smallskip
\Proof Because the quantisation map is linear we have
${\bf c}({\bar\Az})=\sum_{i\ne 1}\;{\bf c}(\Az_{[i]})$. Thus, we
compute
$$\eqalign{{\bf c}({\bar\Az}){\bf c}({\bar\Az})\; =&\;{\bf c}(
\Az^2_{[0]}) +\cr
\noalign{\vskip 0.1cm}   
            &\;[{\bf c}\bigl(\Az_{[0]}\bigr),{\bf c}\bigl(\Az_{[2]}
\bigr)] +\cr
           \noalign{\vskip 0.1cm}
           &\;[{\bf c}\bigl(\Az_{[0]}\bigr),{\bf c}\bigl(\Az_{[3]} 
\bigr)] + \cr
           \noalign{\vskip 0.1cm}
           &\;[{\bf c}\bigl(\Az_{[0]}\bigr),{\bf c}\bigl(\Az_{[4]} 
\bigr)] + {\bf c}(\Az_{[2]})^2 +\cr
           \noalign{\vskip 0.1cm}
           &\ \ \vdots \cr\noalign{\vskip 0.1cm}
           &\; [{\bf c}\bigl(\Az_{[0]}\bigr), {\bf c}\bigl(
\Az_{[i]}\bigr)]+[{\bf c}\bigl(\Az_{[2]}\bigr),{\bf c}\bigl(
\Az_{[i-2]}\bigr)]+\;
                         \dots \cr\noalign{\vskip 0.1cm}
                          &\ \ \vdots\cr}\eqno(2.3)$$
where ${[\;\cdot\; ,\;\cdot\;]}$\ denotes the supercommutator
in $\Omega^0(M,\End\;{\cal E})$. Because  on
the zero-level  ${\bf c}(\Az_{[0]})=
\Az_{[0]}$\ holds together with $\Az_{[0]}\in \Gamma(\End_{C(M)}\;
{\cal E})$\ we get
${\bf c}\bigl(\Az_{[0]}\Az_{[i]}\bigr)=\Az_{[0]}{\bf c}(\Az_{[i]})$\
for $i\ne 1$.
In turn this implies
$[{\bf c}\bigl(\Az_{[0]}\bigr), {\bf c}\bigl(
\Az_{[i]}\bigr)]={\bf c}\bigl([\Az_{[0]}, \Az_{[i]}]\bigr)$\ for all
$i\ne 1$. Finally, by using ${\bf c}({\bar\Az}^2)=
\sum_{i\ne 1}\;{\bf c}\bigl(({\bar\Az}^2)_{[i]}\bigr)$\ together
with (2.2) we obtain the desired result.
\QED 
\smallskip\rm
Thus, this difference ${\bigl({\bf c}({\bar\Az})^2
-{\bf c}({\bar\Az}^2)\bigr)}$\ is independent of the
zero-degree part $\Az_{[0]}$\ of the Clifford superconnection
$\Az$. If $n$\ denotes the
dimension of the underlying manifold $M$, obviously it is true
that $({\bar\Az}^2)_{[i]}=0$\ for $i>n$\ whereas in general
${\bf c}(
\Az_{[i]}){\bf c}(\Az_{[j]})\ne 0$\ for $(i+j)>n$. Consequently
in low dimensions it is  
$\sum_{i,j\ge 2}\;{\bf c}(\Az_{[i]})
{\bf c}(\Az_{[j]})$\ which mainly contributes to the above examined
difference. For example, we have
$$\eqalign{\bigl({\bf c}({\bar\Az})^2
-{\bf c}({\bar\Az}^2)\bigr)=\sum_{i=2}^4
{\bf c} &(\Az_{[i]})^2  -
{\bf c}(\Az_{[2]}^2)+ [{\bf c}\bigl(\Az_{[2]}\bigr),{\bf c}\bigl(
\Az_{[3]}\bigr)] +\cr
& +[{\bf c}\bigl(\Az_{[2]}\bigr),{\bf c}\bigl(
\Az_{[4]}\bigr)] + 
[{\bf c}\bigl(\Az_{[3]}\bigr),{\bf c}\bigl(
\Az_{[4]}\bigr)] \cr}\eqno(2.4)$$
for a Clifford superconnection
$\Az:=\sum_{i=0}^4\;\Az_{[i]}$\ on a Clifford module $\cal E$\
over a four-dimensional Riemannian manifold $M$.
\vfill\break
{\mittel 3. Some canonical constructions}
\ukneu
\vskip 0.7cm
Given a Clifford module $\cal E$\ over a Riemannian mannifold $M$\
we establish in this section a canonical projection
${\bf g}\colon \Omega^*(M,\End\;{\cal E})\rightarrow \Omega^1(M,\End\;
{\cal E})$. This map enables us to
relate a Clifford superconnection $\Az\colon \Gamma({\cal E})
\rightarrow \Omega^*(M,{\cal E})$\
with an ordinary connection
$\widetilde{\nabla}^{\cal E}$\ which in general is not even a Clifford
connection.
However it turns out that this attribution $\Az\mapsto \widetilde{
\nabla}^{\cal E}$\ preserves the most important information
of $\Az$\ concerning our purpose\fussnote{${^2)}$}{
In fact we will show in the next section that
the Dirac operators
$D_\az:={\bf c}\circ \Az$\ and
$\widetilde{D}:= c\circ \widetilde{\nabla}^{\cal E}$\
defined by the Clifford superconnection and the
the corresponding connection, respectively, coincide.}.    
For explaining this we first observe the
\Lemma Let $\nabla^{\cal E}\colon
\Gamma({\cal E})\rightarrow \Gamma(T^*M\otimes {\cal E})$\ be a
Clifford connection on the Clifford module ${\cal E}$.
Then there exists a covariant constant
section $\sigma\in \Gamma(T^*M\otimes \End {\cal E})$\ such that
$c(\sigma)=\eins_{\cal E}$.
\smallskip
\Proof Let $Sym^2(T^*M)$\ denote the bundle of symmetric two
tensors of $T^*M$\ over $M$. There are natural inclusions
${Sym^2(T^*M)\;\hookr{i}\; T^*M\otimes C(M)\;\hookr{j} \;T^*M\otimes
\End {\cal E}}$\
which are compatible with the induced connections,
i.e.
$${j_*\bigl(i_*(\nabla_X s)\bigr)=j_*\bigl(\nabla^{T^*M\otimes C(M)}_X
i_*(s)\bigr)= \nabla^{T^*M\otimes {\petit\rm End}{\cal E}}
j_*\bigl(i_*(s)\bigr)
}\eqno(3.1)$$
for all sections $s\in \Gamma(Sym^2(T^*M))$\ and all
$X\in \Gamma(TM)$. Here the covariant derivative
$\nabla\colon \Gamma(Sym^2(T^*M))
\rightarrow \Gamma(T^*M\otimes Sym^2(T^*M))$\ is induced by the
Levi-Civita connection on $T^*M$\ and the map $i_*\colon
\Gamma(Sym^2(T^*M))\rightarrow
\Gamma(T^*M\otimes C(M))$\
or $j_*\colon \Gamma(T^*M\otimes C(M))\rightarrow \Gamma(T^*M\otimes
\End {\cal E})$\ denote push-forward, respectively.
The second identity holds because we have $\nabla^{T^*M\otimes {\petit
\rm End}{\cal E}}:=\nabla\otimes \eins_{\cal E} + \eins_{T^*M}
\otimes \nabla^{{\petit
\rm End}{\cal E}}$\
with $\nabla^{{\petit\rm End}{\cal E}}$\ induced by the given
Clifford connection $\nabla^{\cal E}$. 
\smallskip
Now take a covariant constant section $\omega\in \Gamma(Sym^2(T^*M))$\  
which yields a non-zero constant $r$\ when evaluated
with the Riemannian metric $g$, i.e $ev_g(\omega)=r$.
Because of (3.1) it is true that $\nabla^{T^*M
\otimes {\petit\rm End}{\cal E}}
j_*(i_*\omega)=0$\ holds. Furthermore we get
$(c\circ j_*\circ i_*)=c^2$\ with
the map $c^2\colon T^*M\otimes T^*M\rightarrow C(M)$\ defined by
$c^2(v\otimes w):=c(v)c(w)$\ for
all $v,w\in \Gamma(T^*M)$. Using the well-known identity
$c^2(v\otimes w)={\bf c}(v\wedge w)-ev_g(v\otimes w)$\ where
${\bf c}\colon \Gamma(\Lambda^*T^*M)
\rightarrow \Gamma(C(M))$\ denotes the quantisation map this can
further be simplified as $(c\circ j_*\circ i_*)=-ev_g$.
So $\sigma:=-{1\over r}\;j_*(i_*\omega)$\ has the
desired properties.
\QED
\smallskip\rm
With respect to a local coordinate system
and using the Einstein convention we may also write
$\sigma=-{1\over r}\;\omega_{\mu\nu}\;dx^\mu\otimes c(dx^\nu)$\
where the coefficients $\omega_{\mu\nu}$\ are totally
symmetric. Note that in general there may exist many
sections of $T^*M\otimes \End {\cal E}$\ with the above
mentioned properties, so by no way we can achieve uniqueness
of the $\End\;{\cal E}$-vallued one-form $\sigma$\ 
in the previous lemma. However there is a canonical choice: If we take
$g\in \Gamma(Sym^2(T^*M))$\ the Riemann metric then, by the above
construction, 
${\gamma:=-{1\over n}\;g_{\mu\nu}\;dx^\mu\otimes c(dx^\nu)\in
\Omega^1(M,\End{\cal E})}$\
with $n:={\rm dim}\;M$\ accomplishes $\nabla^{T^*M\otimes
{\petit\rm End}{\cal E}}\gamma=0$\ and $c(\gamma)=\eins_{\cal E}$.  
In addition, $\gamma$\ is canonical with respect to the Riemannian
structure of $M$.
\smallskip
For any element $\alpha\in \Omega^l(M,\End {\cal E})$\ let
$\mu(\alpha)\colon \Omega^k(M,\End {\cal E})
\rightarrow \Omega^{k+l}(M,\End {\cal E})$\ denote
multiplication in the graded algebra $\Omega^*(M,\End {\cal E})$.
So
the canonical one-form $\gamma\in \Omega^1(M,\End {\cal E})^+$\
induces an even map $\mu(\gamma)\colon\Omega^*(M,\End {\cal E})^\pm
\rightarrow \Omega^{*+1}(M,\End {\cal E})^\pm$.  
On the zero level $\mu(\gamma)\colon \Omega^0(M,\End {\cal E})
\rightarrow \Omega^1(M,\End {\cal E})$\ is injective because
Lemma 3.1 implies that $c\circ\mu(\gamma)=\eins_{\cal E}$\ holds.
We now define 
$${{\bf g}\colon \Omega^*(M,\End {\cal E})\;\mapright{\bf c}\;
\Omega^0(M,\End {\cal E})\;\mapright{\mu(\gamma)}\;
\Omega^1(M,\End {\cal E}).}\eqno(3.2)$$
This map is linear and satisfies
$c\circ {\bf g}=c\circ \mu(\gamma)\circ {\bf c}= \eins_{\cal E}
\circ {\bf c}={\bf c}$\
by the above mentioned property of $\mu(\gamma)$.
Consequently we have ${\bf g}^2={\bf g}
\circ {\bf g}= \mu(\gamma)\circ c \circ {\bf g}
=\mu(\gamma)\circ {\bf c}={\bf g}$, so {\bf g} is a projection.
Let `$\cdot$' be the pointwise defined product in the
algebra bundle $T(M)\otimes \End\;{\cal E}$,
where $T(M)$\ denotes the tensor bundle of $T^*M$\ 
and $i_X\colon \Omega^*(M,\End\;
{\cal E})\rightarrow \Omega^{*-1}(M,\End\;{\cal E})$\ be the inner
derivative with respect to $X\in \Gamma(TM)$.
We will now study in greater detail
this map {\bf g}. Moreover, the following
two lemmas are essential in in proving our theorem.
\Lemma Let $\cal E$\ be a Clifford module over an even-dimensional
Riemannian manifold $M$, 
$\nabla^{\cal E}$\ be a Clifford connection,
$d^{\nabla^{End\;{\cal E}}}\colon
\Omega^*(M,\End\;{\cal E})\rightarrow \Omega^{* + 1}(M,\End\;
{\cal E})$\ be the induced exterior covariant derivative and
$\alpha\in \Omega^*(M,\End\;{\cal E})^-$. Then  
the canonical map ${\bf g}\colon \Omega^*(M,\End\;{\cal E})
\rightarrow \Omega^1(M,\End\;{\cal E})$\
has the properties
$$\eqalignno{{\bf c}\bigl({\bf g}(\alpha)^2\bigr) +
ev_g\bigl({\bf g}(\alpha)\cdot {\bf g}(\alpha)\bigr)
&= {\bf c}
(\alpha)^2 + 2 ev_g\bigl(\beta(\alpha)\cdot
{\bf g}(\alpha)\bigr)
& (3.3)\cr
c^2\bigl(\nabla^{T^*M\otimes {\petit\rm End}{\cal E}}{\bf g}(\alpha)
\bigr) &= {\bf c}\bigl(d^{\nabla^{End{\cal E}}}\alpha\bigr)
-ev_g\nabla^{T^*M\otimes {\petit\rm End}{\cal E}}\beta(\alpha)
& (3.4) \cr}$$
where $\beta(\alpha)\in \Omega^1(M,\End\;{\cal E})$\ is defined
by $\beta(\alpha):=dx^k\otimes {\bf c}\bigl(i(\partial_k)\alpha\bigr)$\
with respect to a local coordinate frame.
\smallskip
\Proof First identity (3.3): With respect to a local
coordinate system we may write ${\bf g}(\alpha)=dx^\mu\otimes
{\bf g}(\alpha)_\mu$. Thus,  ${\bf c}\bigl(
{\bf g}(\alpha)^2\bigr)$\ equals ${1\over 4}\;
[c(dx^\mu),c(dx^\nu)][{\bf g}(
\alpha)_\mu, {\bf g}(\alpha)_\nu]$.
Furthermore using
${1\over 4}\;[c(dx^\mu),c(dx^\nu)][\omega_\mu,\omega_\nu]=
c(dx^\mu)\omega_\mu c(dx^\nu)\omega_\nu +g^{\mu\nu}\omega_\mu
\omega_\nu - c(dx^\mu)[\omega_\nu, c(dx^\nu)]\omega_\nu$\ 
which is true for any $\omega=dx^\mu\otimes \omega_\mu\in
\Omega^1(M,\End\;{\cal E})$\ we obtain 
$$\eqalign{{\bf c}\bigl({\bf g}(\alpha)^2\bigr)=
c(dx^\mu){\bf g}(\alpha)_\mu c(dx^\nu)&{\bf g}(\alpha)_\nu +g^{\mu\nu}
{\bf g}(\alpha)_\mu
{\bf g}(\alpha)_\nu \cr
&- c(dx^\mu)[{\bf g}(\alpha)_\mu, c(dx^\nu)]
{\bf g}(\alpha)_\nu.\cr}\eqno(3.5)$$
Because $c(dx^\mu){\bf g}(\alpha)_\mu=(c\circ {\bf g})(\alpha)=
{\bf c}(\alpha)$\ it remains to show that the sum of the last
two terms in (3.5) equals
$2 ev_g\bigl(\beta(\alpha)\cdot
{\bf g}(\alpha)\bigr)- ev_g\bigl({\bf g}(\alpha)\cdot {\bf g}(
\alpha)\bigr)$. This is a consequence of the following lemma 3.3 .
\smallskip
Second identity (3.4): Using the definition of the map {\bf g}
and the compatibility condition
$\nabla^{T^*M\otimes
{\petit\rm End}{\cal E}}_X\mu(\omega_1)\omega_0=
\mu\bigl(\nabla^{T^*M\otimes
{\petit\rm End}{\cal E}}_X\omega_1\bigr)\omega_0 +
\mu(\omega_1)\bigl(\nabla^{{\petit\rm End}{\cal E}}_X\omega_0\bigr)$\
for all $\omega_i\in \Omega^i(M,\End\;{\cal E}),\ i=0,1$\ and
$X\in \Gamma(TM)$\ we compute
$$\eqalign{\nabla^{T^*M\otimes
{\petit\rm End}{\cal E}}{\bf g}(\alpha) &=dx^\mu\otimes
\nabla^{T^*M\otimes
{\petit\rm End}{\cal E}}_\mu(\mu(\gamma){\bf c}(\alpha))\cr
&=\underbrace{dx^\mu\otimes \mu\bigl(\nabla^{T^*M\otimes
{\petit\rm End}{\cal E}}_\mu \gamma\bigr)}_{\nabla^{T^*M\otimes
{\petit\rm End}{\cal E}}\gamma=0}{\bf c}(\alpha)
+ dx^\mu\otimes \mu(\gamma)\nabla^{{\petit\rm End}{\cal E}}_\mu
{\bf c}(\alpha)\cr}\eqno(3.6)$$
Thus, we get $c^2\bigl(\nabla^{T^*M\otimes
{\petit\rm End}{\cal E}}{\bf g}(\alpha)
\bigr)=c(dx^\mu)c(\mu(\gamma))\nabla^{{\petit\rm End}{\cal E}}_\mu
{\bf c}(\alpha)= c(dx^\mu)\nabla^{{\petit\rm End}{\cal E}}_\mu
{\bf c}(\alpha)$\ by the above mentioned
property of $\mu(\gamma)$.
Note, that because of being induced by a Clifford connection,
$\nabla^{{\petit\rm End}{\cal E}}$\ is compatible with the quantisation
map {\bf c}. More precisely it is true that
$\nabla^{{\petit\rm End}{\cal E}}_X {\bf c}(\omega)={\bf c}\bigl(
\nabla^{\Lambda^*(T^*M)\otimes {\petit\rm End}{\cal E}}_X\omega
\bigr)$\ holds for all forms $\omega\in \Omega^*(M,\End {\cal E})$\ and
$X\in \Gamma(TM)$. Here the tensor product connection
$\nabla^{\Lambda^*(T^*M)\otimes {\petit\rm End}{\cal E}}$\ is defined
by
$\nabla^{\Lambda^*(T^*M)\otimes {\petit\rm End}{\cal E}}:=
\nabla^{\Lambda^*(T^*M)}\otimes \eins_{\cal E} + \eins_{\Lambda^*(T^*M)}
\otimes \nabla^{\petit\rm End{\cal E}}$\ with
the connetion $\nabla^{\Lambda^*(T^*M)}\colon
\Gamma(\Lambda^*(T^*M))\rightarrow
\Gamma(T^*M\otimes \Lambda^*(T^*M))$\ on the exterior bundle being
induced by the Levi-Civita
connection. 
This compatibility condition together
with the equivariance of the quantisation map
${\bf c}\colon \Lambda^*(T^*M)\rightarrow C(M)$\ with
respect to the respective Clifford actions enables us to transform
$$\eqalign{c(dx^\mu)\nabla^{{\petit\rm End}{\cal E}}_\mu
{\bf c}(\alpha) &= c(dx^\mu) {\bf c}\bigl(
\nabla^{\Lambda^*(T^*M)\otimes {\petit\rm End}{\cal E}}_\mu\alpha
\bigr)\cr
&={\bf c}\bigl(dx^\mu\wedge\nabla^{\Lambda^*(T^*M)
\otimes {\petit\rm End}{\cal E}}_\mu\alpha\bigr) -
{\bf c}\bigl(g^{\mu\sigma}i(\partial_\sigma)\nabla^{\Lambda^*(T^*M)
\otimes {\petit\rm End}{\cal E}}_\mu\alpha\bigr)\cr}$$
In the first term by definition
$dx^\mu\wedge\nabla^{\Lambda^*(T^*M)
\otimes {\petit\rm End}{\cal E}}_\mu=:d^{\nabla^{End{\cal E}}}$\ is the
exterior covariant derivative. Thus, it only remains to look after
the second one: Obviously the well-known identity
$[\nabla_\mu, i(\partial_\sigma)]=i(\nabla_\mu\partial_\sigma)$\
implies
$${[\nabla^{\Lambda^*(T^*M)\otimes {\petit\rm End}
{\cal E}}_\mu, i(\partial_\sigma)]=i(\nabla_\mu\partial_\sigma)
=\Gamma^\nu_{\mu\sigma} i(\partial_\nu)}\eqno(3.7)$$
where $\Gamma^\nu_{\mu\sigma}$\ denotes the Cristoffel symbols.
Consequently we obtain
$$\eqalign{{\bf c}\bigl(g^{\mu\sigma}i(\partial_\sigma)\nabla^{\Lambda^*(T^*M)
\otimes {\petit\rm End}{\cal E}}_\mu\alpha\bigr) &=
g^{\mu\sigma}\nabla^{{\petit\rm End}{\cal E}}_\mu
{\bf c}\bigl(i(\partial_\sigma)\alpha\bigr)-g^{\mu\sigma}
\Gamma^\nu_{\mu\sigma} {\bf c}\bigl(i(\partial_\nu)\alpha\bigr)\cr
&=g^{\mu\sigma}\nabla^{T^*M\otimes {\petit\rm End}{\cal E}}_\mu
{\bf c}\bigl(i(\partial_\sigma)\alpha\bigr)}$$
by using the compatibility of the connection $\nabla^{{\petit\rm End}
{\cal E}}$\ with the quantisation map {\bf c} in the inverse
direction.
\QED 
\Lemma Let $\alpha\in \Omega^*(M,\End\;{\cal E})^-$\ and
$\beta(\alpha)\in \Omega^1(M,\End\;{\cal E})^-$\ be defined as
in the previous lemma. Then
${\beta(\alpha)={\bf g}(\alpha)-{1\over 2}\;g_{\sigma\nu}dx^\sigma
\otimes c(dx^\mu)[{\bf g}(\alpha)_\mu, c(dx^\nu)]}$\
holds with respect to a local coordinate frame.
\smallskip
\Proof Using the
property $c\circ {\bf g}={\bf c}$\ of the canonical
projection map {\bf g} we obtain
$${{\bf g}(\alpha)-{1\over 2}\;g_{\sigma\nu}dx^\sigma
\otimes c(dx^\mu)[{\bf g}(\alpha)_\mu, c(dx^\nu)]=dx^\sigma\otimes
\Bigl(-{1\over 2}\;g_{\sigma\nu}[c(dx^\nu), {\bf c}( \alpha)]\Bigr).
}\eqno(3.8)$$ 
Note that $[\;\cdot\; , \;\cdot\;]$\ denotes the supercomutator
in $\Omega^0(M,\End\;{\cal E})$. Thus, it remains to show
$-{1\over 2}\;g_{\sigma\nu}[c(dx^\nu), {\bf c}(\alpha)] =
{\bf c}\bigl(i(\partial_\sigma)\alpha\bigr)$. Because $\End\;{\cal E}
\cong C(M){\hat\otimes} \End_{C(M)}{\cal E}$\ this identity
in $\End\;{\cal E}$\ ensues from the following diagramm
$$\matrix{\Lambda^*(T^*M)\; &\mapright{i(X)}\; &\Lambda^*(T^*M)\cr
\noalign{\vskip 0.2cm}
\mapdown{\bf c} &\quad &\mapdown{\bf c} \cr
\noalign{\vskip 0.2cm}
C(M) & \mapright{-{1\over 2}\;[c(X^*),\;\cdot\;]} &C(M)
\cr}\eqno(3.9)$$
which is commutative for all $X\in \Gamma(TM)$. Here
$X^*:=g(X,\;\cdot\;)\in \Gamma(T^*M)$\ denotes the dual of the
vectorfield $X$\ and $[\;\cdot\; ,\;\cdot\;]$\ is the supercommutator
in the Clifford bundle $C(M)$.
\QED 
\smallskip\rm
Thus, $-c(dx^\mu)[{\bf g}(\alpha)_\mu,c(dx^\nu)]{\bf g}(
\alpha)_\nu=2 g^{\sigma\nu}{\bf c}\bigl(i(\partial_\sigma)\alpha
\bigr){\bf g}(\alpha)_\nu - 2g^{\sigma\nu}{\bf g}(\alpha)_\sigma
{\bf g}(\alpha)_\nu$\ holds and completes the proof of equation (3.3) in
Lemma 3.2 .  
\smallskip
In the following we denote by ${\cal C}({\cal E})$\
the collection of all connections and by
${\cal CSC}({\cal E})$\ the collection of all Clifford
superconnections on a Clifford module $\cal E$. It is well-known
that  
${\cal C}({\cal E})$\ and ${\cal CSC}({\cal E})$\ are affine spaces
modelled over $\Omega^1(M,\End\;{\cal E})$\ and
$\Omega^*(M,\End_{C(M)}{\cal E})$, respectively. With the 
map ${\bf g}\colon \Omega^*(M,\End\;{\cal E})\rightarrow
\Omega^1(M,\End\;{\cal E})$\ in hand we are able to define an affine map
$$\eqalign{{\cal CSC}({\cal E})\;&\longrightarrow\;{\cal C}({\cal E})\cr
\Az &\longmapsto\;\nabla^\az:=\nabla^{\cal E}+{\bf g}({\bar\Az}).
\cr}\eqno(3.10)$$
As before (cf. section 2) we use the
notation $\nabla^{\cal E}:=\Az_{[1]}$\ for the
connection part and ${\bar\Az}$\ for the connection-free
part of $\Az$. Now we compare the corresponding curvature $R^{
\nabla^\az}$\ with $\Fz(\Az)$: By definition, resp. lemma 2.1 we
have
$$\eqalignno{R^{\nabla^\az} &\; =\;R^{\nabla^{\cal E}}\;+\;
d^{\nabla^{End\;E}} {\bf g}({\bar\Az})\; +\; {\bf g}({\bar\Az})^2
&(3.11)\cr
\Fz(\Az) &\;=\;R^{\nabla^{\cal E}}\; +\;
d^{\nabla^{End\;E}} {\bar \Az}\; +\; ({\bar\Az})^2, &(3.12)\cr}$$
and therefore $\bigl(R^{\nabla^\az}-\Fz(\Az)\bigr)= 
(d^{\nabla^{End\;E}} {\bf g}({\bar\Az}) -
d^{\nabla^{End\;E}} {\bar\Az})
+ ({\bf g}({\bar\Az})^2- {\bar\Az}^2)$. Thus, using
lemma 3.2 we obtain the
\Corollary Let $\Az$\ be a Clifford superconnection on a Clifford module
${\cal E}$\ and $\nabla^\az$\ the associated connection with respect to
the affine map (3.10). Then
$$\eqalign{{\bf c}\bigl(R^{\nabla^\az}-\Fz(\Az)\bigr)=
\bigl( {\bf c}( &{\bar \Az} )^2-{\bf c}({\bar \Az}^2)
\bigr)  +ev_g\Bigl(
\nabla^{T^*M\otimes {\petit\rm End}{\cal E}} \bigl({\bf g}({\bar\Az}) -
\beta({\bar\Az})\bigr)\Bigr) +\cr
& + 2ev_g\bigl(\beta({\bar\Az})\cdot {\bf g}({\bar\Az})\bigr)
- ev_g\bigl({\bf g}({\bar \Az})\cdot
{\bf g}({\bar \Az})\bigr)\cr}$$
holds where 
$\Fz(\Az):=\Az^2$\ denotes the supercurvature and
$R^{\nabla^\az}:=(\nabla^\az)^2$\ the curvature of $\Az$\ and $\nabla^\az$,
respectively. 
\rm
\vskip 1cm
{\mittel 4. The generalized Lichnerowicz formula}
\ukneu
\vskip 0.7cm
Let ${\cal E}$\ be a Clifford module over an
even-dimensonal Riemannian manifold $M$.
Generalizing Dirac's original notion a Dirac operator
acting on sections of ${\cal E}$\ can be defined as an odd-parity first order
differential operator $D\colon \Gamma({\cal E}^\pm)\rightarrow
\Gamma({\cal E}^\mp)$\ such that its square $D^2$\ is a generalized
laplacian (cf. [5]). We will regard only those Dirac operators
$D$\ that are compatible with the given Clifford module structure on
${\cal E}$, i.e. $[D, f]=c(df)$\ holds for all $f\in C^\infty(M)$. Note
that this property fully characterizes those Dirac operators. Given
any superconnection $\Az\colon \Gamma({\cal E})^\pm\rightarrow
\Omega^*(M,{\cal E})^\mp$\ on ${\cal E}$, the first order operator
$D_\az$\ defined by the following composition
$${\Gamma({\cal E}^\pm)\;\mapright{\az}\;\Omega^*(M,{\cal E})^\mp\;
\mapright{\cong}\;\Gamma(C(M)\otimes {\cal E})^\mp\;\mapright{c}\;
\Gamma({\cal E}^\mp)}\eqno(4.1)$$
obviously is a Dirac operator. Here the isomorphism is induced by
the quantisation map ${\bf c}\colon \Lambda^*(T^*M)\mapright{\cong} C(M)$\
and the last map denotes the given Clifford action of $C(M)$\ on
${\cal E}$. Note that any connection $\widetilde{\nabla}^{\cal E}\colon
\Gamma({\cal E}^\pm)\rightarrow \Gamma^(T^*M\otimes {\cal E}^\pm)$\ which
respects the grading is also a superconnection. 
Moreover, as it is shown in [5], due to
the above construction (4.1) any Dirac operator is uniquely determined
by a Clifford superconnection, i.e. the assignment
$\Az\mapsto D_\az$\  for $\Az\in {\cal CSC}({\cal E})$\ is a bijection. 
\smallskip
Going back to the map (3.10) we can associate a connection $\nabla^\az\in {
\cal C}({\cal E})$\ to any
Clifford superconnection $\Az\in {\cal CSC}({\cal E})$\ and we
are interessted in comparing the corresponding Dirac operators:     
\Lemma Let $\Az$\ be a Clifford superconnection on a Clifford module
${\cal E}$\ and $\nabla^\az$\ the associated connection with respect to
the affine map (3.10). Then the corresponding Dirac operators
$D_\az$\ and $D_{\nabla^\az}$\ coincide.
\smallskip
\Proof Because we have defined
$\nabla^\az:= \nabla^{\cal E}+ {\bf g}({\bar\Az})$\ where
$\nabla^{\cal E}:=\Az_{[1]}$\ denotes the connection part of the
Clifford superconnection $\Az$\ we obviously obtain
$${D_{\nabla^\az}:=c\circ\nabla^\az=c\circ\Az_{[1]} +(c\circ {\bf g})
({\bar\Az})= c\circ\Az_{[1]}+ {\bf c}({\bar\Az})={\bf c}(\Az)=:D_\az.}$$
Here we have used the property $c\circ {\bf g}={\bf c}$\ of the canonical
projection {\bf g}, cf. the previous section 3.
\QED
Thus, for any Dirac operator $D_\Az$\ on a Clifford module
${\cal E}$\ there exists a connection
$\widetilde{\nabla}^{\cal E}\colon \Gamma({\cal E}^\pm)\rightarrow
\Gamma(T^*M\otimes {\cal E}^\pm)$\
such that
$D_{\widetilde{\nabla}^{\cal E}}=D_\az$. This is a restatement of
Quillen's principle that `Dirac operators are a quantisation of the
theory of connections' (cf. the introduction of [5]). Note that for
the same reason as we remark after lemma 3.1 we can not achieve
uniqueness of the connection $\widetilde{\nabla}^{\cal E}$. However,
the above defined connection $\nabla^\az$\ is the canonical choice. 
\smallskip
If $\widetilde{\nabla}^{\cal E}\colon
\Gamma({\cal E}^\pm)\rightarrow \Gamma(T^*M\otimes
{\cal E}^\pm)$\ is a connection on a Clifford module
${\cal E}$,
recently it has been shown that
there is the following decomposition formula
for the square of the corresponding Dirac operator
${D_{\widetilde{\nabla}^{\cal E}}:= c\circ
\widetilde{\nabla}^{\cal E}}$\ (cf. [2]):
$${D^2_{\widetilde{\nabla}^{\cal E}}
=\triangle^{{\hat\nabla}^{\cal E}}+
{\bf c}(R^{\widetilde{\nabla}^{\cal
E}})+
ev_g\widetilde{\nabla}^{T^*M\otimes {\petit\rm End}{\cal E}}
\varpi_{\widetilde{\nabla}^{\cal E}} +
ev_g(\varpi_{\widetilde{\nabla}^{\cal E}}
\cdot \varpi_{\widetilde{\nabla}^{\cal E}}).}\eqno(4.2)$$
Here $\varpi_{\widetilde{\nabla}^{\cal E}}:=-{1\over 2} g_{\nu\kappa}
dx^\nu\otimes
c(dx^\mu)\bigl([\widetilde{\nabla}^{\cal E}_\mu, c(dx^\kappa)] +
c(dx^\sigma)\Gamma^\kappa_{\sigma\mu}\bigr)\in \Omega^1(M,\End {\cal E})$\
indicates the deviation of the connection
$\widetilde{\nabla}^{\cal E}$\ being a Clifford connection and
$\triangle^{{\hat\nabla}^{\cal E}}$\ is the connection laplacian
associated to 
${\hat\nabla}^{\cal E}:=\widetilde{\nabla}^{\cal E}+
\varpi_{\widetilde{\nabla}^{\cal E}}$\fussnote{${^3)}$}{
With respect
to a local coordinate frame of $TM$, the connection laplacian
${\triangle^{{\hat\nabla}^{\cal E}}}$\
is explicitly
given by ${\triangle^{{\hat\nabla}^{\cal E}}=-g^{\mu\nu}(
{\hat\nabla}^{\cal E}_\mu
{\hat\nabla}^{\cal E}_\nu -\Gamma^\sigma_{\mu\nu}
{\hat\nabla}^{\cal E}_\sigma)}$.}. Only the second term which denotes
the image of the curvature $R^{\widetilde{\nabla}^{\cal E}}
\in \Omega^2(M,
\End({\cal E}))$\ of the given connection $\widetilde{\nabla}^{\cal E}$\
under
the quantisation map
${{\bf c}\colon \Lambda^*T^*M\rightarrow C(M)}$, is endowed with
geometric significance. Of course, if
$\widetilde{\nabla}^{\cal E}$\
is a Clifford connection,
obviously $\varpi_{\widetilde{\nabla}^{\cal E}}=0$\ and therefore
(4.2) reduces to Lichnerowicz's formula
${D_{\widetilde{\nabla}^{\cal E}}^2=\triangle^{
\widetilde{\nabla}^{\cal E}}+{r_M\over 4} +
{\bf c}(R^{{\cal E}/S}_{\widetilde{\nabla}^{\cal E}})}$. Because
any Clifford superconnection $\Az\in {\cal CSC}({\cal E})$\
uniquely determines a Dirac operator $D_{\az}$\ as already
mentioned above, it is natural to reformulate the generalized
Lichnerowicz formula (4.2):
\Theorem Let $\Az=\Az_{[1]}+{\bar\Az}$\ be a Clifford
superconnection on a
Clifford module ${\cal E}$\ over an even-dimensional Riemannian
manifold $M$\ and let $D_\az:={\bf c}\circ \Az$\ denote the
corresponding Dirac operator. Then
$${D_\az^2=\triangle^{{\hat\nabla}^{\cal E}}+{r_M\over 4}+
{\bf c}\bigl(\Fz(\Az)^{{\cal E}/S}\bigr) +
{\bf c}\bigl({\bar\Az}\bigr)^2 - {\bf c}\bigl(
{\bar\Az}^2\bigr) + ev_g\bigl(\beta({\bar\Az})\cdot
\beta({\bar\Az})\bigr)}\eqno(4.3)$$
where 
${\hat\nabla}^{\cal E}:=\Az_{[1]}+\beta({\bar\Az})$\   
determines the connection laplacian $\triangle^{{\hat
\nabla}^{\cal E}}$, $\Fz(\Az)^{{\cal E}/S}$\ denotes the
twisting supercurvature of $\Az\in {\cal CSC}({\cal E})$\ and
$\beta({\bar\Az})\in \Omega^1(M,\End\;
{\cal E})$\ is defined by $\beta({\bar\Az}):= dx^k\otimes
{\bf c}\bigl(i(\partial_k){\bar\Az}\bigr)$\ with respect to a
local coordinate frame.
\smallskip
\Proof For convinience, let $\nabla^{\cal E}:=\Az_{[1]}$\
denote the Clifford connection
part of the Clifford superconnection $\Az$.
Lemma 4.1 tells us that the Dirac operator $D_\az$\
corresponding to $\Az$\ can be equivalently obtained by
$D_\Az=c\circ \nabla^\az$\ using the associated connection
$\nabla^\az:=\nabla^{\cal E}+
{\bf g}({\bar\Az})$. Thus, we reformulate
the decomposition formula (4.2) (cf. [2]):
$${D_\az^2=\triangle^{{\hat\nabla}^{\cal E}}+
{\bf c}(R^{\nabla^A})+
ev_g{\bar\nabla}^{T^*M\otimes {\petit\rm End}{\cal E}}
\varpi_{\nabla^A} +
ev_g(\varpi_{\nabla^A}
\cdot \varpi_{\nabla^A}).}\eqno(4.4)$$
Here we have ${\bar\nabla}^{T^*M\otimes
{\petit\rm End}\;{\cal E}}:=
\nabla\otimes\eins_{\cal E}+\eins_{T^*M}\otimes \widetilde{\nabla}^\az$\
where $\nabla$\ denotes the Levi-Civita connection on $T^*M$\ and
the connection $\widetilde{\nabla}^\az\colon \Gamma(\End{\cal E})
\rightarrow \Gamma(T^*M\otimes\End{\cal E})$\ on the
endomorphism bundle is induced by
$\nabla^\az$.  
\smallskip
Now we inspect the right hand side of formula (4.4) term by term:
Recall that 
$\nabla^{\cal E}$\ is compatible with the Clifford action, i.e.
$[\nabla^{\cal E}_\mu,c(dx^\kappa)]=-c(dx^\sigma)
\Gamma^\kappa_{\sigma\mu}$\ holds with respect to a local
coordinate frame. This implies
$$\eqalign{\varpi_{\nabla^A}:&=-{1\over 2}\;g_{\nu\kappa}
dx^\nu\otimes c(dx^\mu)[(\nabla^\az_\mu-\nabla^{\cal E}_\mu),c(dx^\kappa)
]\cr
&=-{1\over 2}\;g_{\nu\kappa}
dx^\nu\otimes c(dx^\mu)[{\bf g}({\bar\Az})_\mu),c(dx^\kappa)]\cr
}\eqno(4.5)$$
and therefore ${\hat\nabla}^{\cal E}:=\nabla^\az+\varpi_{\nabla^A}=
\nabla^{\cal E}+\bigl({\bf g}({\bar\Az})-{1\over 2}\;g_{\nu\kappa}
dx^\nu \otimes c(dx^\mu)[{\bf g}({\bar\Az})_\mu),c(dx^\kappa)]\bigr)$\
holds. Thus, using lemma 3.3 it is true that
${\hat\nabla}^{\cal E}=\Az_{[1]}+\beta({\bar\Az})$\ holds with 
$\beta({\bar\Az})\in \Omega^1(M,\End\;
{\cal E})$\ locally defined by $\beta({\bar\Az}):= dx^k\otimes
{\bf c}\bigl(i(\partial_k){\bar\Az}\bigr)$.
\smallskip
Before we replace the curvature term ${\bf c}\bigl(R^{\nabla^A}
\bigr)$\ by
${\bf c}\bigl(\Fz(\Az)\bigr)+ {\bf c}\Bigl(R^{\nabla^A}-\Fz(\Az)\bigr)$\
which involves the supercurvature, we study the last two terms in
(4.4). Note that, again by lemma 3.3, we have $\varpi_{\nabla^A}=
\beta({\bar\Az})-{\bf g}({\bar\Az})$. Applying
${\bar\nabla}^{T^*M\otimes {\petit\rm End}\;{\cal E}}_\mu=\nabla^{T^*M
\otimes {\petit \rm End}\;{\cal E}}_\mu+[{\bf g}(
{\bar\Az})_\mu,\;\cdot\;]$\
where $[\;\cdot\; ,\;\cdot\;]$\
denotes the commutator in $\End\;{\cal E}$\ and the tensor
connection $\nabla^{T^*M\otimes {\petit\rm End}\;{\cal E}}:=
\nabla\otimes \eins_{\cal E}+ \eins_{T^*M}\otimes \nabla^{\cal E}$\
is defined as in the proof of lemma 3.1 we obtain for
the third term   
$$\eqalign{ev_g{\bar\nabla}^{T^*M\otimes {\petit\rm End}{\cal E}}
\varpi_{\nabla^A}=ev_g&\Bigl(\nabla^{T^*M
\otimes {\petit \rm End}\;{\cal E}}
\bigl(\beta({\bar\Az})-{\bf g}({\bar\Az})\bigr)\Bigr)+\cr
&+ ev_g\bigl({\bf g}({\bar\Az})\cdot \beta({\bar\Az})\bigr)-ev_g\bigl(
\beta({\bar\Az})\cdot {\bf g}({\bar\Az})\bigr).}\eqno(4.6)$$
For the forth one we calculate
$$\eqalign{ev_g(\varpi_{\nabla^A}\cdot \varpi_{\nabla^A})&=ev_g\Bigl(
\bigl(\beta({\bar\Az})-{\bf g}({\bar\Az})\bigr)\cdot
\bigl(\beta({\bar\Az})-{\bf g}({\bar\Az})\bigr)\Bigr)\cr
& = ev_g\bigl(\beta({\bar\Az})\cdot \beta({\bar\Az})\bigr)-
ev_g\bigl({\bf g}({\bar\Az})\cdot \beta({\bar\Az})\bigr)-\cr
&\ \ \ \ -ev_g\bigl(\beta({\bar\Az})\cdot {\bf g}({\bar\Az})\bigr)+
ev_g\bigl({\bf g}({\bar\Az})\cdot {\bf g}({\bar\Az})\bigr),\cr
}\eqno(4.7)$$
always provided that the dot `$\cdot$' indicates the fibrewise defined
product in the algebra bundle $T(M)\otimes \End\;{\cal E}$\ with
$T(M)$\ being the tensor bundle of $T^*M$. Clearly,
$$\eqalign{\hbox{\rm (4.6) + (4.7)}=ev_g &\Bigl(\nabla^{T^*M
\otimes {\petit \rm End}\;{\cal E}}
\bigl(\beta({\bar\Az})-{\bf g}({\bar\Az})\bigr)\Bigr)
-2 ev_g\bigl(\beta({\bar\Az})\cdot {\bf g}({\bar\Az})\bigr)+\cr
&+ev_g\bigl({\bf g}({\bar\Az})\cdot {\bf g}({\bar\Az})\bigr)+
ev_g\bigl(\beta({\bar\Az})\cdot \beta({\bar\Az})\bigr)\cr}\eqno(4.8)$$
If we additionally insert
${{\bf c}\bigl(\Fz(\Az)\bigr)+ {\bf c}\Bigl(R^{\nabla^A}-\Fz(\Az)\Bigr)}$\
for ${\bf c}\bigl(R^{\nabla^A}\bigr)$\ in formula (4.4) and use
corollary 3.4, the first three terms of (4.8) cancel out.
Finally, by using proposition 3.43 of [BGV] the supercurvature $\Fz(\Az)$\
decomposes under the
isomorphism $\End\;{\cal E}\cong C(M)\otimes\End\;{\cal E}$\
as $\Fz(\Az)=c(R)+\Fz(\Az)^{{\cal E}/S}$\ where
$c(R)\in \Omega^2(M,C(M))$\ is the action of the Riemannian
curvature $R$\ of $M$\ on the Clifford module and $\Fz(\Az)^{{\cal E}/S}
\in \Omega^*(M,\End_{C(M)}{\cal E})$\ denotes the twisting
supercurvature of $\Az$. This completes the proof of the
theorem.
\QED
Obviously, the
endomorphism $P({\bar\Az})\in \Gamma(\End\;{\cal E})$\ defined by
the last three terms in the decomposition formula 
$P({\bar\Az}):={\bf c}\bigl({\bar\Az}\bigr)^2 - {\bf c}\bigl(
{\bar\Az}^2\bigr) + ev_g\bigl(\beta({\bar\Az})\cdot
\beta({\bar\Az})\bigr)$\ 
depends only on the higher degree parts 
$\Az_{[i]},\ i\ge 2$\ of the Clifford superconnection
$\Az$. This means $P({\bar\Az})=0$\ for $\Az=\Az_{[0]}+\Az_{[1]}$.
Hence the generalized Lichnerowicz formula
reduces to $D_\az^2=\triangle^{\nabla^{\cal E}}+{r_M\over 4}+
{\bf c}\bigl(\Fz(\Az)^{{\cal E}/S}\bigr)$\
in this case. Furthermore, if we denote
$\Az_{[0]}:=\Phi\in \End_{C(M)}^-{\cal E}$, 
we obtain
$${D_{(\Phi+\nabla^{\cal E})}^2=
\triangle^{\nabla^{\cal E}}
+{r_M\over 4} +
{\bf c}(R^{{\cal E}/S}_{\nabla^{\cal E}}) +
c\;\nabla^{{\petit\rm End}{\cal E}}(\Phi) +
\Phi^2.}\eqno(4.9)$$
Thus, we recover the decomposition formula for the square of
a `Dirac operator of simple type', cf. [2].
\smallskip   
As we have already mentioned in the introduction, Getzler 
has stated the generalized Lichnerowicz formula in
the form $D_\az^2=\triangle^{{\hat\nabla}^{\cal E}}+{r_M\over 4}+
{\bf c}\bigl(\Fz(\Az)^{{\cal E}/S}\bigr) + P(\Az)$. In contrast
to our result (4.3) above, 
he has not specified the endomorphism $P(\Az)\in \Gamma(\End\;{\cal E})$\
in general. However, calculation of $P(\Az)$\ for 
$\Az:=\Az_{[0]}+\Az_{[1]}+\Az_{[2]}$\ where the
two-form part is given by
$\Az_{[2]}:={1\over 2}\;dx^i\wedge dx^j\otimes \omega_{ij}$\
with $\omega_{ij}\in \End_{C(M)}{\cal E}$\ for all $i,\;j$\
yields
$${P(\Az)=2g^{ij}\;{\bf c}(dx^k\wedge dx^l)\;\omega_{ik}\omega_{jl}
- g^{ij} g^{kl}\; \omega_{ik}\omega_{jl}}\eqno(4.10)$$
which reproduces Getzler's example (cf. [7]).     
\vskip 1cm
$\hbox{\mittel 5. Families of Dirac operators and the extended}$\par
$\hskip 0.6cm\hbox{\mittel   generalized Lichnerowicz formula}$
\ukneu
\vskip 0.7cm  
Suppose that $\pi\colon M\rightarrow B$\ is a family of oriented
even-dimensional
Riemannian manifolds $(M_z\;\vert\; z\in B)$\ and ${\cal E}$\
is a bundle over $M$\ such that its fibrewise restriction
${\cal E}_z:={\cal E}\vert_z$\ is a Clifford module for each $z\in B$.
Furthermore, let $\nabla^{\cal E}$\ be a connection on ${\cal E}$\
which is a Clifford connection when restricted to each Clifford module
${\cal E}_z$\ and let $\hbox{\sq D}:=(\hbox{\sq D}^z\;\vert\; z\in B)$\
be the associated family of Dirac operators.  
In [4], Bismut constructed a superconnection $\nAz$\ 
on the $C^\infty$-direct image $\pi_*{\cal E}$\ corresponding to the
family $\hbox{\sq D}$\ whose Chern character, by using
Getzler's rescaling trick, is explicitly computable. 
A crucial step in this calculation was the remarkably simple
formula for the supercurvature $\nAz^2$\ (cf. [4] and [5], chapter
10.3) which can be understood as extending the ordinary
Lichnerowicz formula to this infinit-dimensional context.  
In this section we follow Bismut's construction with $\nabla^{\cal E}$\
replaced
by an arbitrary superconnection $\Az$\ on $\cal E$\ whose
restriction to each bundle ${\cal E}_z$\ is a Clifford superconnection.
This enables us to introduce a superconnection $\nAz^\az$\
on $\pi_*{\cal E}$\ corresponding to the family
of Dirac operators $\hbox{\sq D}_\az:=
(\hbox{\sq D}^z_\az\;\vert\;
z\in B)$\ defined by $\Az$. For the supercurvature
$(\nAz^\az)^2$\ we also obtain a simple formula. Similarly, that
extends the above generalisation (4.3)
of Lichnerowicz's decomposition. 
\smallskip
First, we recall briefly the geometric structure (cf. [5]):
Let $\pi\colon M\rightarrow B$\ as above with a metric
$g_{M/B}$\ on the vertical tangent bundle $T(M/B)$\
and a connection on $TM$. This
induces a decomposition $TM=T(M/B)\oplus \pi^*TB$\
where we have identified $\pi^*TB$\ with the horizontal space.
Let $\nabla^{M/B}$\ denote the connection on $T(M/B)$\ associated
to the vertical metric $g_{M/B}$\ as constructed in [4].
Given a metric $g_B$\ on the base $B$\ with associated
Levi-Civita connection $\nabla^B$\ on $TB$\ 
one is able to define $\nabla^{\oplus}:=\nabla^{M/B}\oplus
\pi^*\nabla^B$. Note that $\nabla^\oplus$\ preserves the
metric $g=\pi^*g_B\oplus g_{M/B}$\ but
differs from
the corresponding Levi-Civita connection $\nabla^g$\
by a torsion-term, i.e. $(\nabla^g-\nabla^\oplus)\in
\Omega^1(M, Sk(TM))$\ where $Sk(TM)$\ denotes the bundle
of skew-symmetric endomorphisms of $TM$. Using the
isomorphism $\tau\colon \Lambda^2 TM \mapright{\cong} Sk(TM)$\ which is
given by the Riemannian metric $g$\ we have $\nabla^g=\nabla^\oplus
+{1\over 2}\;\tau(\omega)$. Interessting, this
form $\omega\in \Omega^1(M,\Lambda^2 TM)$\
is independent of the choosen metric $g_B$, cf. [4].
Thus, using the blow-up metric $g_u:=u^{-1} \pi^*g_B+ g_{M/B}$\ with
$u\in \rz\setminus \{0\}$, the limit as $u\rightarrow 0$\ of
the corresponding Levi-Civita
connections $\nabla^{M,u}$\ exits and yields
$${\nabla^{M,0}\;:=\;\lim_{u\rightarrow 0}\;\nabla^{M,u}\;=\;
\nabla^g,}\eqno(5.1)$$
cf. [B] or [BGV].
Now let $g^u=u\pi^*g_B\oplus g_{M/B}$\ denote the
corresponding dual metric
and $\nabla^{T^*M,u},\nabla^{T^*M,0}$\ the corresponding
dual connections on $T^*M$. Note that the limit $g^0:=
\lim_{u\rightarrow 0}\;g^u$\ is degenerate on the horizontal
space $\pi^*T^*B$. Furthermore, suppose that ${\cal E}$\ is a hermitian vector bundle
over $M$\ which, in addition, is a Clifford module along the
fibres of the
bundle $(M,B,\pi)$. More precisely, there is a Clifford action
$c\colon C(M/B)\rightarrow \End\;{\cal E}$\ where $C(M/B)$\ denotes
the bundle of Clifford algebras over $M$\ generated by the vertical
bundle $(T^*(M/B), g_{M/B})$. We now define a natural action $m_0$\ of
$C_0(M):=\pi^*\Lambda T^*B\otimes C(M/B)$\ on the bundle
$\Ee:=\pi^*\Lambda T^*B\otimes {\cal E}$\ over $M$:
$${m_0(a):=\cases{\epsilon(a)\otimes id_{\cal E} &iff $a\in
\pi^* T^*B$ \cr
id_{\pi^*\Lambda T^*B}\otimes c(a) &iff $a\in T^*(M/B)$.\cr}}\eqno(5.2)$$
Here $\epsilon$\ denotes exterior multiplication on
$\pi^*\Lambda T^*B$.
Note that this action
$m_0$\ can be understood as the limit
$\lim_{u\rightarrow 0}\;m_u$\ of the Clifford actions
$m_u$\ of the Clifford bundles
$C_u(M)$\ generated by $T^*M$\ with respect to the relation
$v\star w + w\star v=-2g^u(v,w)$\ for all $v,w \in \Gamma(T^*M)$\ on
$\Ee$\ as defined in [5]. Thus, in this reference $m_0$\
is called `degenerate
Clifford action'.
\smallskip
Now let $\Az\colon \Gamma({\cal E}^\pm)\rightarrow \Omega^*(M,
{\cal E})^\mp$\ be a superconnection which is a Clifford
superconnection with respect to the above defined Clifford action
$c$\ of the vertical Clifford bundle $C(M/B)$, i.e.
$[\Az,c(a)]=c(\nabla^{M/B}a)$\ holds for all $a\in C(M/B)$. Then
$$\eqalignno{\Az^{\Ee,\oplus} &:=\pi^*\nabla^B\otimes \eins_{\cal E}
+\eins_{\pi^*\Lambda T^*B}\otimes \Az &(5.3)\cr
\Az^{\Ee,u} &:= \Az^{\Ee, \oplus} + {1\over 2}\;m_u(\omega) &(5.4)\cr
}$$
are superconnections on
$\Ee$\ where $\omega\in \Omega^1(M, \Lambda^2T^*M)$\ denotes
the above mentioned torsion term considered as operating on
$\Ee$\ by the Clifford action $m_u$\fussnote{$^4)$}{Recall the Lie-algebra
isomorphism $\Lambda^2T^*M\cong C_u^2(M):=\{ {\bf m}_u(a)\;\vert\;
a\in \Lambda^2T^*M\}$.}. Here ${\bf m}_u\colon \Lambda T^*M\rightarrow
C_u(M)$\ denotes the respective quantisation maps. Obviously, by a
similar argument as in Proposition 10.10 of [BGV],
the various superconnections $\Az^{\Ee, u}$\ are Clifford superconnections
on $\Ee$\ with
respect to the corresponding Clifford actions $m_u$, i.e.
$[\Az^{\Ee, u}, m_u(a)]=m_u(\nabla^{T^*M,u}a)$\ for all
$a\in C_u(M)$\ holds. Let $D_{\Az,u}:={\bf m}_{\petit\boldit u}\circ
\Az^{\Ee,u}$\ denote the corresponding Dirac operators. Furthermore, taking
the limit  
$${\Az^{\Ee,0}\;:=\;\lim_{u\rightarrow 0}\; \Az^{\Ee, u}\;
 =\; \Az^{\Ee,\oplus}
+ {1\over 2}\;m_0(\omega)}\eqno(5.5)$$
we obtain a superconnection with the property
$[\Az^{\Ee,0}, m_0(a)]=m_0(\nabla^{T^*M,0}a)$\ for all
$a\in C_0(M)$. In other words, $\Az^{\Ee,0}$\ respects the
degenerate Clifford action
$m_0$.
\smallskip
Recall that given a vector bundle ${\cal E}$\ over $M$, the
$C^\infty$-direct image $\pi_*{\cal E}$\ is the infinite-dimensional
bundle over $B$\ whose fibre at $z\in B$\ is defined to be
the space $\Gamma({\cal E}_z)$\ of all $C^\infty$-sections of the bundle
${\cal E}_z$\ over $M_z$. Imitating the key idea in Bismut's construction
we use the isomorphism
$${\Omega^*(B,\pi_*{\cal E})\;\cong \;\Gamma(\pi^*\Lambda T^*B\otimes
{\cal E})}\eqno(5.6)$$
to define an operator $\nAz^\az\colon \Omega^*(B,\pi_*{\cal E})^\pm
\rightarrow \Omega^*(B,\pi_*{\cal E})^\mp$\ by
$${\nAz^\az\; :=\;\lim_{u\rightarrow 0}\;D_{\Az,u}
\;=\;\lim_{u\rightarrow 0}\bigl({\bf m}_{\petit\boldit u}\circ\Az^{\Ee,u}
\bigr)}\eqno(5.7)$$ 
Equivallently we may write $\nAz^\az={\bf m}_{\petit\boldit 0}\circ
\Az^{\Ee,0}$\ if we use definition (5.5). Hence 
$\nAz^\az$\ can be understood as a kind of `Dirac operator' on the bundle
$\Ee$\ provided with the `Clifford module'-structure
which is defined by the degenerate Clifford action $m_0$.
\smallskip
Before studying
this operator
further, recall that  
$S\in \Gamma(T^*_HM\otimes
\End\;T(M/B))\cong \Gamma(T^*_HM\otimes T^*(M/B)\otimes T(M/B)$\
defined by
$S(Z,\theta,X):=\theta(\nabla^{M/B}_ZX -P[Z,X])$\ for
$Z\in \Gamma(T_HM),\; \theta\in \Gamma(T^*(M/B)$\ and
$X\in \Gamma(T(M/B))$\ is the second fundamental form
associated to a family $\pi\colon M\rightarrow B$\
of Riemannian manifolds with a given
splitting $TM=T(M/B)\oplus T_HM$\ and a
connection $\nabla^{M/B}$\ on the
vertical bundle $T(M/B)$, cf. [5]. Here $P\colon
TM\rightarrow T(M/B)$\ denotes the projection map with
kernel the choosen horizontal space. Now we obtain the following
\Lemma The operator
$\nAz^\az\colon \Gamma(\pi_*{\cal E})^\pm
\rightarrow \Omega^*(B,\pi_*{\cal E})^\mp$\ defined by $\nAz^\az:=
{\bf m}_{\petit\boldit 0}\circ \Az^{\Ee,0}$\ is a superconnection on
the direct image $\pi_*{\cal E}$\ which is explicitly given
by
$${\nAz^\az={\bf c}\circ \Az + \epsilon\circ (\Az+ {1\over 2}k)
+{\bf m}_{\petit\boldit 0}\circ \Omega}\eqno(5.8)$$
where $k\in \Omega^1(M)$\ defined
by $k(Z):=tr\bigl(S(Z)\bigr)$\ denotes the mean curvature and
$\Omega\in \Gamma(\Lambda^2T_H^*M\otimes T(M/B))$\ is the curvature
of the connection $\nabla^{M/B}$\ associated to the
family $\pi\colon M\rightarrow B$\ of Riemannian manifolds.
\smallskip
\Proof We follow the proof of Proposition 10.15 of [5]: Obviously,
the operator $\nAz^\az$\ satisfies
$\nAz^\az(\nu s)=(\epsilon\circ \nabla^B\nu)s + (-1)^{\vert\;\nu\;\vert}\nu\nAz^\az s$\
for all $\nu\in \Omega^*(B)$\ and $s\in \Gamma(\Ee)$. Since
$\nabla^B$\ is the
Levi-Civita associated to the choosen metric $g_B$\ on the base and
therefore torsion-free, we see that $\epsilon\circ \nabla^B=d_B$\ is
the exterior covariant derivative on the base. Thus, $\nAz^\az$\
is a superconnection on $\pi_*{\cal E}$.
\smallskip  
For proving the
explicit formula (5.8), we observe that
the splitting of the cotangent bundle $T^*M=\pi^*T^*B\oplus T^*(M/B)$\ implies
$\Lambda T^*M =\pi^*\Lambda T^*B\otimes \Lambda T^*(M/B)$. Thus
the quantisation map
${\bf m}_{\petit\boldit 0}\colon \Lambda T^*M\mapright{\cong} C_0(M):=
\pi^*\Lambda T^*B\otimes C(M/B)$\ associated to the degenerate
Clifford structure on $\Ee$\ is given by
${\bf m}_{\petit\boldit 0}=\epsilon\otimes \eins_{C(M/B)}+
\eins_{\pi^*\Lambda T^*B}\otimes {\bf c}$\ when restricted to
$T^*M$. Furthermore, because of definition (5.5) we know
$${\nAz^\az={\bf c}\circ\Az^{\Ee,\oplus} +\epsilon\circ \Az^{\Ee,\oplus}
+{1\over 2}\;m_0(\omega).}\eqno(5.9)$$
Using Lemma 10.13 of [5] which tells us
$m_0(\omega)= {1\over 2}(\epsilon\circ k)
+{\bf m}_{\petit\boldit 0}\circ \Omega$\ completes the proof
of equation (5.8).
\QED
\smallskip
Now assume that the superconnection $\Az$\ on 
${\cal E}$\ consists only of the connection part $\Az_{[1]}=
\nabla^{\cal E}$. Then, by construction, the corresponding superconnection
$\nAz^{\nabla^{\cal E}}$\
on $\pi_*{\cal E}$\ is the
Bismut superconnection. Hence we call the
operator $\nAz^\Az$\ on $\pi_*{\cal E}$\
defined by (5.7) `the generalized Bismut
superconnection'. In any case $\nAz^\az_{[0]}=
\hbox{\sq D}_\az$\ holds with $\hbox{\sq D}_\az$\
being a family of Dirac operators
defined by the superconnection
$\Az$\ as above. Thus, the generalized
Bismut superconnection $\nAz^\az$\ corresponds to an
arbitrary family of Dirac operators.
\smallskip
Given any connection $\nabla\colon \Gamma(\Ee)\rightarrow \Gamma(
T^*M\otimes \Ee)$\ we will call the second order operator
$\triangle_{M/B}^{\nabla}$\ on $\Ee$\ defined by
$\triangle_{M/B}^\nabla:=g_{M/B}^{ij}\bigl(\nabla_i\nabla_j-
\nabla_{\nabla^{M/B}_i e_j}\bigr)$\ the `vertical' connection
laplacian associated to $\nabla$. Note that we adopt the
convention that the indices $i,\; j,\;\dots $\
label vertical vectors. Using this
notation, we finally state the 
analogue of Theorem 4.2 which provides the formula for 
supercurvature 
$\bigl(\nAz^\az\bigr)^2$\ of the generalized Bismut
superconnection:
\Theorem Let $\nAz^\az\colon \Gamma(\pi_*{\cal E})^\pm\rightarrow
\Omega^*(B,\pi_*{\cal E})$\ be the generalized Bismut superconnection
corresponding to the Clifford superconnection
$\Az$\ on the Clifford module ${\cal E}$\ over $C(M/B)$. Then
$${\!\bigl(\nAz^\az\bigr)^2\!=\triangle_{\!M/B\!}^{\hat\nabla}\! +\!
 {r_{\!M/B\!}
\over 4}
+{\bf m}_{\petit\boldit 0}\!\bigl(\Fz(\Az)^{{\cal E}/S}\bigr) +
{\bf m}_{\petit\boldit 0}\!\bigl({\bar\Az}\bigr)^2\!-\!
{\bf m}_{\petit\boldit 0}\!\bigl({\bar\Az}^2\bigr) + ev_{g_{M/B} }\!
\bigl(\beta_0({\bar\Az})\cdot\beta_0({\bar \Az})\bigr)}$$
where $r_{M/B}$\ denotes the scalar curvature of $\nabla^{M/B}$,
${\hat\nabla}:=\Az_{[1]}+\beta_0({\bar\Az})$\ determines the
`vertical' connection laplacian $\triangle_{M/B}^{\hat\nabla}$, the
form
$\Fz(\Az)^{{\cal E}/S}\in \Omega^*(M,\End_{C(M/B)}{\cal E})$\ is
the twisting curvature  of $\Az$\ and   
$\beta_0({\bar\Az})\in \Omega^1(M,\End\;\Ee)$\ is
defined by $\beta_0({\bar\Az}):= dx^j\otimes
{\bf m}_{\petit\boldit 0}\bigl(i(\partial_j){\bar\Az}\bigr)$\
with respect to a
local vertical coordinate frame.
\smallskip\rm
The proof of this theorem is similar to our proof of
the generalized Lichnerowicz formula and can be found in [1]. However,
in order to do so, the extension of all
the previous results in section 3 to this case is indispensable. 
For instance, the analogue ${\bf g}_{\petit\bf M/B}$\
of the canonical projection
map ${\bf g}$, cf. (3.2), can be defined by
$${{\bf g}_{\petit\bf M/B}\colon \Omega^*(M,\End \;\Ee)\;
\mapright{{\petit\bf m}_{\bf 0}}\;\Omega^0(M,\End\;\Ee)\;
\mapright{\mu(\gamma_{M/B})}\;\Omega^1(M,\End\;\Ee).}\eqno(5.10)$$   
Recall that here ${\bf m}_{\petit\bf 0}$\ denotes the
composition of the `quantisation map' ${\bf m}_{\petit\bf 0}$\ and the
action $m_0\colon \pi^*\Lambda T^*B{\hat\otimes} C(M/B)
\rightarrow \End\;\Ee$\ by abuse of notation, and the
one-form
$\gamma_{M/B}$\ is
defined by $\gamma_{M/B}=-{1\over n} (g_{M/B})_{ij} dx^i\otimes c(dx^j)\in
\Omega^1(M,\End\;\Ee)$\ where $n:= {\rm dim}\;M_z,\ z\in B$\
is the fibre dimension. The $\End\;\Ee$-vallued one forms decompose
into components
$\Omega^1(M,\End\;\Ee)\cong
\Gamma(\pi^*T^*B\otimes \End\;\Ee)\oplus \Gamma(T^*(M/B)\otimes \End\;
\Ee)$\ and we calculate
$m_0\circ {\bf g}_{\petit\bf M/B}=(\epsilon\otimes\eins_\Ee +
c\otimes \eins_\Ee)\circ \mu(\gamma_{M/B})\circ {\bf m}_{\petit\bf
0}=c\circ \mu(\gamma_{M/B})\circ {\bf m}_{\petit\bf 0}=
{\bf m}_{\petit\bf 0}$\ because $\gamma_{M/B}$\ is
vertical\fussnote{${^5)}$}{Recall, that $c\colon C(M/B)\rightarrow
\End\;{\cal E}$\ denotes the vertical Clifford action.}. 
In turn this implies an analogue of Lemma 4.1 above.
Thus, generalized Bismut superconnections
can be understood as  `quantisation' (with respect to the
`degenerate' quantisation map ${\bf m}_{\petit\boldit 0}$)
of the theory of connections on the `degenerate Clifford module'
$\Ee$. This is an extension of Quillen's principle as
mentioned after Lemma 4.1 
concerning the relation of Dirac operators and connections to
the setting of families of
Dirac operators. For more details and a complete proof of
Theorem 5.2 we refer once more to [1] where it is also shown how
to compute the Chern character of a generalized Bismut superconnection
using this formula.

\rm
\vskip 1cm
{\mittel 6. Conclusion}
\ukneu
\vskip 0.7cm
We have studied Dirac operators acting on sections of
a Clifford module ${\cal E}$\ over a Riemannian manifold $M$. Motivated
by the fact that any Clifford superconnection $\Az$\ on ${\cal E}$\
uniquely determines a Dirac operator $D_\az$, in this paper
we have emphasized the supersymmeric approach using Quillen's
super-formalism. We have proven the supersymmetric version (4.3)
of the decomposition formula for the square of a Dirac operator $D_\Az$ 
which generalizes the classical result [8] due to Lichnerowicz.
Associated to a family of (arbitrary) Dirac operators
$\hbox{\sq D}_\az:=\{D_\az^z\;\vert\; z\in B\}$\
parametrized by a  not necessarily finite
dimensional manifold $B$\ we have defined the
notion of `generalized Bismut superconnection' $\nAz^\az$. This
generalizes Bismut's construction [4]. Similarly we have obtained
a simple formula for its supercurvature $\bigl(\nAz^\az\bigr)^2$,
extending the generalized Lichnerowicz formula (4.3). This might be
seen as a first step to prove the local Atiyah-Singer index
theorem also for families of arbitrary Dirac operators [1]. For
applications of the generalized Lichnerowicz formula in physics,
we refer to [2] and [3].
\vskip 0.7cm
{\bf Acknowledgements.} Many thanks to Barbara
and to my parents for all their understanding, love and support. 
I also want to express my gratitude to E. Binz for the possibility
of further using the equipment of his professorial chair
and his gentle encouragement.
\vskip 1cm
\begref
\ref{[1]} T.Ackermann, \sl The local family index theorem revisited\rm ,
to appear\smallskip
\ref{[2]} T.Ackermann, J.Tolksdorf, \sl The generalized
Lichnerowicz formula and analysis of Dirac operators, to
appear in Journ. f\"ur die reine u. angewandte Mathem.\smallskip 
\ref{[3]} T.Ackermann, J.Tolksdorf, \sl A Unification of
Gravity and Yang-Mills-Higgs Gauge Theories, to appear in Letters of
Math. Phys. \smallskip
\ref{[4]} J.M. Bismut, \sl The Atiyah-Singer
index theorem for families of Dirac operators: Two heat equation
proofs, \rm Inv. Mathematicae \bf 83\rm (1989), 92\smallskip
\ref{[5]} N.Berline, E.Getzler, M.Vergne, \sl Heat kernels and
Dirac operators\rm , Springer (1992)\smallskip
\ref{[6]} E.Getzler, \sl A short proof of the local Atiyah-Singer
index theorem\rm , Topology {\bf 25} (1986), 111-117\smallskip
\ref{[7]} E.Getzler, \sl The Bargmann Representation,
Generalized Dirac Operators and the Index of
Pseudodifferential Operators on $\rz^n$\rm , Cont. Mathematics
\bf 179\rm (1994), 63-81\smallskip
\ref{[8]} A.Lichnerowicz, \sl Spineurs harmonique\rm ,
C. R. Acad. Sci. Paris S${\acute {\rm e}}$r. A {\bf 257} (1963)
\smallskip
\ref{[9]} D.Quillen, \sl Superconnections and the Chern Character\rm ,
Toplology {\bf 24} (1985), 89-95\vfill
\bye